

\input amstex
\expandafter\ifx\csname mathdefs.tex\endcsname\relax
  \expandafter\gdef\csname mathdefs.tex\endcsname{}
\else \message{Hey!  Apparently you were trying to
  \string twice.   This does not make sense.} 
\errmessage{Please edit your file (probably \jobname.tex) and remove
any duplicate ``\string\input'' lines} \fi




\catcode`\X=12\catcode`\@=11

\def\n@wcount{\alloc@0\count\countdef\insc@unt}
\def\n@wwrite{\alloc@7\write\chardef\sixt@@n}
\def\n@wread{\alloc@6\read\chardef\sixt@@n}
\def\r@s@t{\relax}\def\v@idline{\par}\def\@mputate#1/{#1}
\def\l@c@l#1X{\firstpart.#1}\def\gl@b@l#1X{#1}\def\t@d@l#1X{{}}

\def\crossrefs#1{\ifx\all#1\let\tr@ce=\all\else\def\tr@ce{#1,}\fi
   \n@wwrite\cit@tionsout\openout\cit@tionsout=\jobname.cit 
   \write\cit@tionsout{\tr@ce}\expandafter\setfl@gs\tr@ce,}
\def\setfl@gs#1,{\def\@{#1}\ifx\@\empty\let\next=\relax
   \else\let\next=\setfl@gs\expandafter\xdef
   \csname#1tr@cetrue\endcsname{}\fi\next}
\def\m@ketag#1#2{\expandafter\n@wcount\csname#2tagno\endcsname
     \csname#2tagno\endcsname=0\let\tail=\all\xdef\all{\tail#2,}
   \ifx#1\l@c@l\let\tail=\r@s@t\xdef\r@s@t{\csname#2tagno\endcsname=0\tail}\fi
   \expandafter\gdef\csname#2cite\endcsname##1{\expandafter
     \ifx\csname#2tag##1\endcsname\relax?\else\csname#2tag##1\endcsname\fi
     \expandafter\ifx\csname#2tr@cetrue\endcsname\relax\else
     \write\cit@tionsout{#2tag ##1 cited on page \folio.}\fi}
   \expandafter\gdef\csname#2page\endcsname##1{\expandafter
     \ifx\csname#2page##1\endcsname\relax?\else\csname#2page##1\endcsname\fi
     \expandafter\ifx\csname#2tr@cetrue\endcsname\relax\else
     \write\cit@tionsout{#2tag ##1 cited on page \folio.}\fi}
   \expandafter\gdef\csname#2tag\endcsname##1{\expandafter
      \ifx\csname#2check##1\endcsname\relax
      \expandafter\xdef\csname#2check##1\endcsname{}%
      \else\immediate\write16{Warning: #2tag ##1 used more than once.}\fi
      \multit@g{#1}{#2}##1/X%
      \write\t@gsout{#2tag ##1 assigned number \csname#2tag##1\endcsname\space
      on page \number\count0.}%
   \csname#2tag##1\endcsname}}
\def\multit@g#1#2#3/#4X{\def\t@mp{#4}\ifx\t@mp\empty%
      \global\advance\csname#2tagno\endcsname by 1 
      \expandafter\xdef\csname#2tag#3\endcsname
      {#1\number\csname#2tagno\endcsnameX}%
   \else\expandafter\ifx\csname#2last#3\endcsname\relax
      \expandafter\n@wcount\csname#2last#3\endcsname
      \global\advance\csname#2tagno\endcsname by 1 
      \expandafter\xdef\csname#2tag#3\endcsname
      {#1\number\csname#2tagno\endcsnameX}
      \write\t@gsout{#2tag #3 assigned number \csname#2tag#3\endcsname\space
      on page \number\count0.}\fi
   \global\advance\csname#2last#3\endcsname by 1
   \def\t@mp{\expandafter\xdef\csname#2tag#3/}%
   \expandafter\t@mp\@mputate#4\endcsname
   {\csname#2tag#3\endcsname\lastpart{\csname#2last#3\endcsname}}\fi}
\def\t@gs#1{\def\all{}\m@ketag#1e\m@ketag#1s\m@ketag\t@d@l p
   \m@ketag\gl@b@l r \n@wread\t@gsin
   \openin\t@gsin=\jobname.tgs \re@der \closein\t@gsin
   \n@wwrite\t@gsout\openout\t@gsout=\jobname.tgs }
\outer\def\localtags{\t@gs\l@c@l}
\outer\def\globaltags{\t@gs\gl@b@l}
\outer\def\newlocaltag#1{\m@ketag\l@c@l{#1}}
\outer\def\newglobaltag#1{\m@ketag\gl@b@l{#1}}

\newif\ifpr@ 
\def\m@kecs #1tag #2 assigned number #3 on page #4.%
   {\expandafter\gdef\csname#1tag#2\endcsname{#3}
   \expandafter\gdef\csname#1page#2\endcsname{#4}
   \ifpr@\expandafter\xdef\csname#1check#2\endcsname{}\fi}
\def\re@der{\ifeof\t@gsin\let\next=\relax\else
   \read\t@gsin to\t@gline\ifx\t@gline\v@idline\else
   \expandafter\m@kecs \t@gline\fi\let \next=\re@der\fi\next}
\def\pretags#1{\pr@true\pret@gs#1,,}
\def\pret@gs#1,{\def\@{#1}\ifx\@\empty\let\n@xtfile=\relax
   \else\let\n@xtfile=\pret@gs \openin\t@gsin=#1.tgs \message{#1} \re@der 
   \closein\t@gsin\fi \n@xtfile}

\newcount\sectno\sectno=0\newcount\subsectno\subsectno=0
\newif\ifultr@local \def\ultralocal{\ultr@localtrue}
\def\firstpart{\number\sectno}
\def\lastpart#1{\ifcase#1 \or a\or b\or c\or d\or e\or f\or g\or h\or 
   i\or k\or l\or m\or n\or o\or p\or q\or r\or s\or t\or u\or v\or w\or 
   x\or y\or z \fi}

\def\resetall{\global\advance\sectno by 1\subsectno=0
   \gdef\firstpart{\number\sectno}\r@s@t}
\def\resetsub{\global\advance\subsectno by 1
   \gdef\firstpart{\number\sectno.\number\subsectno}\r@s@t}
\def\newsection#1\par{\resetall\vskip0pt plus.3\vsize\penalty-250
   \vskip0pt plus-.3\vsize\bigskip\bigskip
   \message{#1}\leftline{\bf#1}\nobreak\bigskip}
\def\subsection#1\par{\ifultr@local\resetsub\fi
   \vskip0pt plus.2\vsize\penalty-250\vskip0pt plus-.2\vsize
   \bigskip\smallskip\message{#1}\leftline{\bf#1}\nobreak\medskip}

\def\t@gsoff#1,{\def\@{#1}\ifx\@\empty\let\next=\relax\else\let\next=\t@gsoff
   \def\@@{p}\ifx\@\@@\else
   \expandafter\gdef\csname#1cite\endcsname##1{\zeigen{##1}}
   \expandafter\gdef\csname#1page\endcsname##1{?}
   \expandafter\gdef\csname#1tag\endcsname##1{\zeigen{##1}}\fi\fi\next}
\def\verbatimtags{\ifx\all\relax\else\expandafter\t@gsoff\all,\fi}
\def\zeigen#1{\hbox{$\langle$}#1\hbox{$\rangle$}}

\def\(#1){\edef\dot@g{\ifmmode\ifinner(\hbox{\noexpand\etag{#1}})
   \else\noexpand\eqno(\hbox{\noexpand\etag{#1}})\fi
   \else(\noexpand\ecite{#1})\fi}\dot@g}

\newif\ifbr@ck
\def\eat#1{}
\def\[#1]{\br@cktrue[\br@cket#1'X]}
\def\br@cket#1'#2X{\def\temp{#2}\ifx\temp\empty\let\next\eat
   \else\let\next\br@cket\fi
   \ifbr@ck\br@ckfalse\br@ck@t#1,X\else\br@cktrue#1\fi\next#2X}
\def\br@ck@t#1,#2X{\def\temp{#2}\ifx\temp\empty\let\neext\eat
   \else\let\neext\br@ck@t\def\temp{,}\fi
   \def\teemp{#1}\ifx\teemp\empty\else\rcite{#1}\fi\temp\neext#2X}
\def\resetbr@cket{\gdef\[##1]{[\rtag{##1}]}}
\def\references{\resetbr@cket\newsection References\par}

\newtoks\symb@ls\newtoks\s@mb@ls\newtoks\p@gelist\n@wcount\ftn@mber
    \ftn@mber=1\newif\ifftn@mbers\ftn@mbersfalse\newif\ifbyp@ge\byp@gefalse
\def\defm@rk{\ifftn@mbers\n@mberm@rk\else\symb@lm@rk\fi}
\def\n@mberm@rk{\xdef\m@rk{{\the\ftn@mber}}%
    \global\advance\ftn@mber by 1 }
\def\rot@te#1{\let\temp=#1\global#1=\expandafter\r@t@te\the\temp,X}
\def\r@t@te#1,#2X{{#2#1}\xdef\m@rk{{#1}}}
\def\b@@st#1{{$^{#1}$}}\def\str@p#1{#1}
\def\symb@lm@rk{\ifbyp@ge\rot@te\p@gelist\ifnum\expandafter\str@p\m@rk=1 
    \s@mb@ls=\symb@ls\fi\write\f@nsout{\number\count0}\fi \rot@te\s@mb@ls}
\def\byp@ge{\byp@getrue\n@wwrite\f@nsin\openin\f@nsin=\jobname.fns 
    \n@wcount\currentp@ge\currentp@ge=0\p@gelist={0}
    \re@dfns\closein\f@nsin\rot@te\p@gelist
    \n@wread\f@nsout\openout\f@nsout=\jobname.fns }
\def\m@kelist#1X#2{{#1,#2}}
\def\re@dfns{\ifeof\f@nsin\let\next=\relax\else\read\f@nsin to \f@nline
    \ifx\f@nline\v@idline\else\let\t@mplist=\p@gelist
    \ifnum\currentp@ge=\f@nline
    \global\p@gelist=\expandafter\m@kelist\the\t@mplistX0
    \else\currentp@ge=\f@nline
    \global\p@gelist=\expandafter\m@kelist\the\t@mplistX1\fi\fi
    \let\next=\re@dfns\fi\next}
\def\symbols#1{\symb@ls={#1}\s@mb@ls=\symb@ls} 
\def\bigsymbol{\textstyle}
\symbols{\bigsymbol\ast,\dagger,\ddagger,\sharp,\flat,\natural,\star}
\def\ftnumbers{\ftn@mberstrue} \def\ftsymbols{\ftn@mbersfalse}
\def\paginal{\byp@ge} \def\resetftnumbers{\ftn@mber=1}
\def\ftnote#1{\defm@rk\expandafter\expandafter\expandafter\footnote
    \expandafter\b@@st\m@rk{#1}}

\long\def\jump#1\endjump{}
\def\ssum{\mathop{\lower .1em\hbox{$\textstyle\Sigma$}}\nolimits}

\def\qed{\nobreak\kern 1em \vrule height .5em width .5em depth 0em}
\def\newneq{\hbox{\rlap{\hbox to 1\wd9{\hss$=$\hss}}\raise .1em 
   \hbox to 1\wd9{\hss$\scriptscriptstyle/$\hss}}}
\def\subsetne{\setbox9 = \hbox{$\subset$}\mathrel{\hbox{\rlap
   {\lower .4em \newneq}\raise .13em \hbox{$\subset$}}}}
\def\supsetne{\setbox9 = \hbox{$\subset$}\mathrel{\hbox{\rlap
   {\lower .4em \newneq}\raise .13em \hbox{$\supset$}}}}

\def\vbar{\mathchoice{\vrule height6.3ptdepth-.5ptwidth.8pt\kern-.8pt}
   {\vrule height6.3ptdepth-.5ptwidth.8pt\kern-.8pt}
   {\vrule height4.1ptdepth-.35ptwidth.6pt\kern-.6pt}
   {\vrule height3.1ptdepth-.25ptwidth.5pt\kern-.5pt}}
\def\f@dge{\mathchoice{}{}{\mkern.5mu}{\mkern.8mu}}
\def\b@c#1#2{{\rm \mkern#2mu\vbar\mkern-#2mu#1}}
\def\b@b#1{{\rm I\mkern-3.5mu #1}}
\def\b@a#1#2{{\rm #1\mkern-#2mu\f@dge #1}}
\def\bb#1{{\count4=`#1 \advance\count4by-64 \ifcase\count4\or\b@a A{11.5}\or
   \b@b B\or\b@c C{5}\or\b@b D\or\b@b E\or\b@b F \or\b@c G{5}\or\b@b H\or
   \b@b I\or\b@c J{3}\or\b@b K\or\b@b L \or\b@b M\or\b@b N\or\b@c O{5} \or
   \b@b P\or\b@c Q{5}\or\b@b R\or\b@a S{8}\or\b@a T{10.5}\or\b@c U{5}\or
   \b@a V{12}\or\b@a W{16.5}\or\b@a X{11}\or\b@a Y{11.7}\or\b@a Z{7.5}\fi}}

\catcode`\X=11 \catcode`\@=12

\expandafter\ifx\csname citeadd.tex\endcsname\relax
\expandafter\gdef\csname citeadd.tex\endcsname{}
\else \message{Hey!  Apparently you were trying to
\string twice.   This does not make sense.} 
\errmessage{Please edit your file (probably \jobname.tex) and remove
any duplicate ``\string\input'' lines} \fi

\sectno=-1   
\localtags
\ifx\shlhetal\undefinedcontrolsequence\let\shlhetal\relax\fi
\NoBlackBoxes
\define\mr{\medskip\roster}
\define\sn{\smallskip\noindent}
\define\mn{\medskip\noindent}
\define\bn{\bigskip\noindent}
\define\ub{\underbar}
\define\wilog{\text{without loss of generality}}
\define\ermn{\endroster\medskip\noindent}

\define\dbcu{\dsize\bigcup}
\define \nl{\newline}
\documentstyle {amsppt}
\topmatter
\title{Borel Whitehead Groups} \endtitle
\author {Saharon Shelah \thanks {\null\newline I would like to thank 
Alice Leonhardt for the beautiful typing. \null\newline
 First Typed - 98/Mar/4 \null\newline
 Latest Revision -  98/Sept/25 \null\newline
 \S1, \S2 done, Fall '89 \null\newline
 Publication no. 402} \endthanks} \endauthor 
\affil{Institute of Mathematics\\
 The Hebrew University\\
 Jerusalem, Israel
 \medskip
 Rutgers University\\
 Mathematics Department\\
 New Brunswick, NJ  USA
 \medskip
 MSRI \\
 Berkeley, CA  USA} \endaffil
\mn

\abstract {We investigate the Whiteheadness of Borel abelian groups
($\aleph_1$-free, \wilog \, as otherwise this is trivial).
We show that CH (and even WCH) implies any such abelian group is free, and
always $\aleph_2$-free.} \endabstract
\endtopmatter
\document  

\expandafter\ifx\csname alice2jlem.tex\endcsname\relax
  \expandafter\gdef\csname alice2jlem.tex\endcsname{}
\else \message{Hey!  Apparently you were trying to
\string  twice.   This does not make sense.}
\errmessage{Please edit your file (probably \jobname.tex) and remove
any duplicate ``\string\input'' lines} \fi

\expandafter\ifx\csname bib4plain.tex\endcsname\relax
  \expandafter\gdef\csname bib4plain.tex\endcsname{}
\else \message{Hey!  Apparently you were trying to \string twice.   This does not make sense.}
\errmessage{Please edit your file (probably \jobname.tex) and remove
any duplicate ``\string\input'' lines} \fi

\def\renewcommand{\newcommand}	       
\edef\cite{\the\catcode`@}%
\catcode`@ = 11
\let\@oldatcatcode = \cite
\chardef\@letter = 11
\chardef\@other = 12
%
%
%
%
\def\@innerdef#1#2{\edef#1{\expandafter\noexpand\csname #2\endcsname}}%
%
%
\@innerdef\@innernewcount{newcount}%
\@innerdef\@innernewdimen{newdimen}%
\@innerdef\@innernewif{newif}%
\@innerdef\@innernewwrite{newwrite}%
%
%
%
\def\@gobble#1{}%
%
%
%
\ifx\inputlineno\@undefined
   \let\@linenumber = \empty 
\else
   \def\@linenumber{\the\inputlineno:\space}%
\fi
%
%
%
\def\@futurenonspacelet#1{\def\cs{#1}%
   \afterassignment\@stepone\let\@nexttoken=
}%
\begingroup 
\def\\{\global\let\@stoken= }%
\\ 
\endgroup
\def\@stepone{\expandafter\futurelet\cs\@steptwo}%
\def\@steptwo{\expandafter\ifx\cs\@stoken\let\@@next=\@stepthree
   \else\let\@@next=\@nexttoken\fi \@@next}%
\def\@stepthree{\afterassignment\@stepone\let\@@next= }%
%
%
%
\def\@getoptionalarg#1{%
   \let\@optionaltemp = #1%
   \let\@optionalnext = \relax
   \@futurenonspacelet\@optionalnext\@bracketcheck
}%
%
%
\def\@bracketcheck{%
   \ifx [\@optionalnext
      \expandafter\@@getoptionalarg
   \else
      \let\@optionalarg = \empty
      \expandafter\@optionaltemp
   \fi
}%
\def\@@getoptionalarg[#1]{%
   \def\@optionalarg{#1}%
   \@optionaltemp
}%
%
%
%
\def\@nnil{\@nil}%
\def\@fornoop#1\@@#2#3{}%
\def\@for#1:=#2\do#3{%
   \edef\@fortmp{#2}%
   \ifx\@fortmp\empty \else
      \expandafter\@forloop#2,\@nil,\@nil\@@#1{#3}%
   \fi
}%
\def\@forloop#1,#2,#3\@@#4#5{\def#4{#1}\ifx #4\@nnil \else
       #5\def#4{#2}\ifx #4\@nnil \else#5\@iforloop #3\@@#4{#5}\fi\fi
}%
\def\@iforloop#1,#2\@@#3#4{\def#3{#1}\ifx #3\@nnil
       \let\@nextwhile=\@fornoop \else
      #4\relax\let\@nextwhile=\@iforloop\fi\@nextwhile#2\@@#3{#4}%
}%
%
%
%
\@innernewif\if@fileexists
\def\@testfileexistence{\@getoptionalarg\@finishtestfileexistence}%
\def\@finishtestfileexistence#1{%
   \begingroup
      \def\extension{#1}%
      \immediate\openin0 =
         \ifx\@optionalarg\empty\jobname\else\@optionalarg\fi
         \ifx\extension\empty \else .#1\fi
         \space
      \ifeof 0
         \global\@fileexistsfalse
      \else
         \global\@fileexiststrue
      \fi
      \immediate\closein0
   \endgroup
}%
%
%
%
%
\def\bibliographystyle#1{%
   \@readauxfile
   \@writeaux{\string\bibstyle{#1}}%
}%
\let\bibstyle = \@gobble
%
%
\let\bblfilebasename = \jobname
\def\bibliography#1{%
   \@readauxfile
   \@writeaux{\string\bibdata{#1}}%
   \@testfileexistence[\bblfilebasename]{bbl}%
   \if@fileexists
      \nobreak
      \@readbblfile
   \fi
}%
\let\bibdata = \@gobble
%
%
\def\nocite#1{%
   \@readauxfile
   \@writeaux{\string\citation{#1}}%
}%
\@innernewif\if@notfirstcitation
%
%
\def\cite{\@getoptionalarg\@cite}%
%
%
\def\@cite#1{%
   \let\@citenotetext = \@optionalarg
   \printcitestart
   \nocite{#1}%
   \@notfirstcitationfalse
   \@for \@citation :=#1\do
   {%
      \expandafter\@onecitation\@citation\@@
   }%
   \ifx\empty\@citenotetext\else
      \printcitenote{\@citenotetext}%
   \fi
   \printcitefinish
}%
\def\@onecitation#1\@@{%
   \if@notfirstcitation
      \printbetweencitations
   \fi
   \expandafter \ifx \csname\@citelabel{#1}\endcsname \relax
      \if@citewarning
         \message{\@linenumber Undefined citation `#1'.}%
      \fi
      \expandafter\gdef\csname\@citelabel{#1}\endcsname{%
\strut
\vadjust{\vskip-\dp\strutbox
\vbox to 0pt{\vss\parindent0cm \leftskip=\hsize 
\advance\leftskip3mm
\advance\hsize 4cm\strut\openup-4pt 
\rightskip 0cm plus 1cm minus 0.5cm ?  #1 ?\strut}}
         {\tt
            \escapechar = -1
            \nobreak\hskip0pt
            \expandafter\string\csname#1\endcsname
            \nobreak\hskip0pt
         }%
      }%
   \fi
   \csname\@citelabel{#1}\endcsname
   \@notfirstcitationtrue
}%
%
%
\def\@citelabel#1{b@#1}%
%
%
\def\@citedef#1#2{\expandafter\gdef\csname\@citelabel{#1}\endcsname{#2}}%
%
%
%
\def\@readbblfile{%
   \ifx\@itemnum\@undefined
      \@innernewcount\@itemnum
   \fi
   \begingroup
      \def\begin##1##2{%
         \setbox0 = \hbox{\biblabelcontents{##2}}%
         \biblabelwidth = \wd0
      }%
      \def\end##1{}
      %
      %
      \@itemnum = 0
      \def\bibitem{\@getoptionalarg\@bibitem}%
      \def\@bibitem{%
         \ifx\@optionalarg\empty
            \expandafter\@numberedbibitem
         \else
            \expandafter\@alphabibitem
         \fi
      }%
      \def\@alphabibitem##1{%
         \expandafter \xdef\csname\@citelabel{##1}\endcsname {\@optionalarg}%
         \ifx\biblabelprecontents\@undefined
            \let\biblabelprecontents = \relax
         \fi
         \ifx\biblabelpostcontents\@undefined
            \let\biblabelpostcontents = \hss
         \fi
         \@finishbibitem{##1}%
      }%
      \def\@numberedbibitem##1{%
         \advance\@itemnum by 1
         \expandafter \xdef\csname\@citelabel{##1}\endcsname{\number\@itemnum}%
         \ifx\biblabelprecontents\@undefined
            \let\biblabelprecontents = \hss
         \fi
         \ifx\biblabelpostcontents\@undefined
            \let\biblabelpostcontents = \relax
         \fi
         \@finishbibitem{##1}%
      }%
      \def\@finishbibitem##1{%
         \biblabelprint{\csname\@citelabel{##1}\endcsname}%
         \@writeaux{\string\@citedef{##1}{\csname\@citelabel{##1}\endcsname}}%
         \ignorespaces
      }%
      %
      %
      \let\em = \bblem
      \let\newblock = \bblnewblock
      \let\sc = \bblsc
      \frenchspacing
      \clubpenalty = 4000 \widowpenalty = 4000
      \tolerance = 10000 \hfuzz = .5pt
      \everypar = {\hangindent = \biblabelwidth
                      \advance\hangindent by \biblabelextraspace}%
      \bblrm
      \parskip = 1.5ex plus .5ex minus .5ex
      \biblabelextraspace = .5em
      \bblhook
      \input \bblfilebasename.bbl
   \endgroup
}%
%
%
\@innernewdimen\biblabelwidth
\@innernewdimen\biblabelextraspace
%
%
%
\def\biblabelprint#1{%
   \noindent
   \hbox to \biblabelwidth{%
      \biblabelprecontents
      \biblabelcontents{#1}%
      \biblabelpostcontents
   }%
   \kern\biblabelextraspace
}%
%
%
%
\def\biblabelcontents#1{{\bblrm [#1]}}%
%
%
\def\bblrm{\rm}%
%
%
\def\bblem{\it}%
%
%
\def\bblsc{\ifx\@scfont\@undefined
              \font\@scfont = cmcsc10
           \fi
           \@scfont
}%
%
%
\def\bblnewblock{\hskip .11em plus .33em minus .07em }%
%
%
\let\bblhook = \empty
%
%
%
\def\printcitestart{[}
\def\printcitefinish{]}
\def\printbetweencitations{, }
\def\printcitenote#1{, #1}
%
%
%
\let\citation = \@gobble
%
%
%
\@innernewcount\@numparams
%
%
\def\newcommand#1{%
   \def\@commandname{#1}%
   \@getoptionalarg\@continuenewcommand
}%
%
%
\def\@continuenewcommand{%
   \@numparams = \ifx\@optionalarg\empty 0\else\@optionalarg \fi \relax
   \@newcommand
}%
%
%
\def\@newcommand#1{%
   \def\@startdef{\expandafter\edef\@commandname}%
   \ifnum\@numparams=0
      \let\@paramdef = \empty
   \else
      \ifnum\@numparams>9
         \errmessage{\the\@numparams\space is too many parameters}%
      \else
         \ifnum\@numparams<0
            \errmessage{\the\@numparams\space is too few parameters}%
         \else
            \edef\@paramdef{%
               \ifcase\@numparams
                  \empty  No arguments.
               \or ####1%
               \or ####1####2%
               \or ####1####2####3%
               \or ####1####2####3####4%
               \or ####1####2####3####4####5%
               \or ####1####2####3####4####5####6%
               \or ####1####2####3####4####5####6####7%
               \or ####1####2####3####4####5####6####7####8%
               \or ####1####2####3####4####5####6####7####8####9%
               \fi
            }%
         \fi
      \fi
   \fi
   \expandafter\@startdef\@paramdef{#1}%
}%
%
%
%
%
\def\@readauxfile{%
   \if@auxfiledone \else 
      \global\@auxfiledonetrue
      \@testfileexistence{aux}%
      \if@fileexists
         \begingroup
            \endlinechar = -1
            \catcode`@ = 11
            \input \jobname.aux
         \endgroup
      \else
         \message{\@undefinedmessage}%
         \global\@citewarningfalse
      \fi
      \immediate\openout\@auxfile = \jobname.aux
   \fi
}%
%
%
\newif\if@auxfiledone
\ifx\noauxfile\@undefined \else \@auxfiledonetrue\fi
%
%
%
%
\@innernewwrite\@auxfile
\def\@writeaux#1{\ifx\noauxfile\@undefined \write\@auxfile{#1}\fi}%
%
%
%
\ifx\@undefinedmessage\@undefined
   \def\@undefinedmessage{No .aux file; I won't give you warnings about
                          undefined citations.}%
\fi
%
%
\@innernewif\if@citewarning
\ifx\noauxfile\@undefined \@citewarningtrue\fi
%
%
%
\catcode`@ = \@oldatcatcode


\def\widestnumber#1#2{}

\def\rm{\fam0 \tenrm}

\def\fakesubhead#1\endsubhead{\bigskip\noindent{\bf#1}\par}


%
%
%

%

\font\textrsfs=rsfs10
\font\scriptrsfs=rsfs7
\font\scriptscriptrsfs=rsfs5

\newfam\rsfsfam
\textfont\rsfsfam=\textrsfs
\scriptfont\rsfsfam=\scriptrsfs
\scriptscriptfont\rsfsfam=\scriptscriptrsfs

\edef\oldcatcodeofat{\the\catcode`\@}
\catcode`\@11

\def\Cal@@#1{\noaccents@ \fam \rsfsfam #1}

\catcode`\@\oldcatcodeofat

\newpage

\head {\S0 Introduction} \endhead  \resetall
\bigskip

\definition{\stag{0.1} Definition}  1) We say that $\bar \psi = \langle
\psi_0,\psi_1 \rangle$ is a code for a Borel abelian group if:
\mr
\item "{$(a)$}"  $\psi_0(\ldots,\ldots)$ codes a Borel equivalence relation
$E = E^{\bar \psi}$ on a subset $B_* = B^{\bar \psi}_*$ of ${}^\omega 2$ so
$[\psi_0(\eta,\eta) \leftrightarrow \eta \in B_*]$ and
$[\psi_0(\eta,\nu) \rightarrow \eta \in B_* \and \nu \in B_*]$, the group will
have a set of elements $B = B^{\bar \psi}_*/E^{\bar \psi}$
\sn
\item "{$(b)$}"   $\psi_1 = \psi_1(x,y,z)$ code a Borel set of triples from
${}^\omega 2$ such that \nl
$\{(x/E^{\bar \psi},y/E^{\bar \psi},z/E^{\bar \psi}):\psi_1(x,y,z)\}$ 
is the graph of a function from $B \times B$ to $B$ such
that $(B,+)$ is an abelian group.
\ermn
2) We say Borel$^+$ if (b) is replaced by:
\mr
\item "{$(b)'$}"  $\psi_1$ codes a Borel function from $B_* \times B_*$ to
$B_*$ which respects $E^{\bar \psi}$, the function is called $+$ and 
$(B,+)$ is
an abelian group (well, we should denote the function which $+$  induces from
$(B_*/E^{\bar \psi}) \times (B_*/E^{\bar \psi})$ into $B_*/E^{\bar \psi}$ by
e.g. $+_{E^{\bar \psi}}$, but are not strict).
\ermn
We let $B^{\bar \psi} = B_{\bar \psi} = (B,+)$ be the group coded by 
$\bar \psi$; abusing notation we may write $B$ for $B_{\bar \psi}$.
\enddefinition
\bn
Clearly \nl 
\stag{0.1A}\ub{Observation}:  The set of codes for Borel
abelian groups is $\Pi^1_2$. \nl
An abelian group $B$ is Borel if it has a Borel code.
\bn
An interesting problem suggested by Dave Marker is the Borel version of 
Whitehead's problem: namely \nl
\ub{\stag{0.2} Question}:  Is every Borel Whitehead group free?

In this paper we will give a partial answer to this question.  We will show
that every Borel Whitehead group is $\aleph_2$-free.  In particular, the
continuum hypothesis implies that every Borel Whitehead group is free.  This
latter result provides a contrast to the author's proof (\cite{Sh:98}) that
it is consistent with CH that there is a Whitehead group of cardinality
$\aleph_1$ which is not free.

We refer the reader to \cite{EM} for the necesary background material on
abelian groups.
\medskip

Suppose $B$ is an $\aleph_1$-free abelian group.  Let $S_0 = \{G \subset
B:|G| = \aleph_0$ and $B/G$ is not $\aleph_1$-free$\}$.  It is well known
that if $B$ is not
$\aleph_2$-free, then $S_0$ is stationary.  We will argue that the converse
is true for Borel abelian groups and the answer is quite absolute.  
Lastly, we deal with weakening Borel to Souslin.
\bn
\stag{0.3}\ub{Question}:  If $B$ is an $\aleph_2$-free Borel abelian group,
what can be the $n$ in the analysis of a nonfree $\aleph_2$-free abelian
subgroup of $B$ from \cite{Sh:161} (or see \cite{EM} or \cite{Sh:523})?

We thank Todd Eisworth for corrections.
\newpage

\head {\S1 On $\aleph_2$-freeness} \endhead  \resetall
\bigskip

\demo{\stag{1.1} Hypothesis}  Let $B$ be an $\aleph_1$-free Borel abelian
group.  Let $\bar \psi$ be a Borel code for $B$.

Let $S_B = S_{\bar \psi} = 
\{K \subseteq B:K \text{ is a countable subgroup and } B/K \text{ is not }
\aleph_1 \text{-free}\}$.
\enddemo
\bigskip

\proclaim{\stag{1.2} Lemma}  1) If $S_B$ is stationary, \ub{then} $B$ is not
$\aleph_2$-free. \nl
2) Moreover, there is an increasing continuous sequence $\langle G_i:
i < \omega_1 \rangle$ of countable subgroups of $B$ such that $G_{i+1}/G_i$ 
is not free for each $i < \omega_1$.
\endproclaim
\bigskip

\remark{Remark}  On such proof in mode theory see \cite[\S2]{Sh:43},
\cite{BKM78} and \cite{Sch85}.
\endremark
\bigskip

\demo{Proof}  We work in a universe $V \models ZFC$.  Force with 
$\bold P = \{p:p \text{ is a function from}$ \nl
some $\alpha < \omega_1 \text{ to } {}^\omega 2\}$.
Let $G \subseteq \bold P$ be $V$-generic and let $V[G]$ denote the generic
extension.

Since $\bold P$ is $\aleph_1$-closed, forcing with $\bold P$ adds no new
reals.  Thus $\bar \psi$ still codes $B$ in the generic extension, i.e.
$B^{V[G]}_{\bar \psi} = B^V_{\bar \psi}$.  
Forcing with $\bold P$ also adds no new countable subsets 
of $B$ hence ``$B$ is $\aleph_1$-free" holds in $V$ iff it holds in $V[G]$.
Similarly if $K \subset B$ is countable, then ``$B/K$ is $\aleph_1$-free" 
holds in $V$ iff it holds in $V[G]$.  Thus, $S^V_{\bar \psi} = 
S^{V[G]}_{\bar \psi}$.  Moreover,  since $\bold P$ is proper,
$S_{\bar \psi}$ remains stationary (see \cite[Ch.III]{Sh:f}).

Since $V[G] \models CH$, we can write

$$
B = \dbcu_{\alpha < \omega_1} B_\alpha,
$$
\mn
where $\bar B = \langle B_\alpha:\alpha < \omega_1 \rangle$ is an increasing
continuous chain of countable subgroups.  Let $S = \{\alpha < \omega_1:
B/B_\alpha$ is not $\aleph_1$-free$\}$.  Since $S_{\bar \psi}$ is 
stationary (as a subset of $[B]^{\aleph_0}$) necessarily, $S$ is a stationary
subset of $\omega_1$.  So $V[G] \models$ ``$B$ is not free".

By Pontryagon's criteria for each $\alpha \in S$ there are $n_\alpha \in
\omega$ and $a^\alpha_0,\dotsc,a^\alpha_{n_\alpha}$ such that

$$
PC(B_\alpha \cup \{a^\alpha_0,\dotsc,a^\alpha_{n_\alpha}\})/B_\alpha
$$
\mn
is not free, 
where $PC(X) = PC(X,B)$ is the pure closure of the subgroup of $B$
which $X$ generates.  We choose $n_\alpha$ minimal with this property.

Work in $V[G]$.  Let $\kappa$ be a regular cardinal such that ${\Cal H}
(\kappa)$ satisfies enough axioms of set theory to handle all of our 
arguments, and let $<^*$ be a well ordering of ${\Cal H}(\kappa)$.  Let $N
\preceq ({\Cal H}(\kappa),\in,<^*)$ be countable such that $\bar \psi,S,
\langle B_\alpha:\alpha < \omega_1 \rangle$ and $\left < \langle 
a^\alpha_0,\dotsc,a^\alpha_{n_\alpha} \rangle:\alpha < \omega_1 \right>$ 
belong to $N$.

The model $N$ has been built in $V[G]$, but since forcing with $\bold P$
adds no new reals, there is a transitive model $N_0 \in V$ isomorphic to $N$
and let $h$ be an isomorphism from $N$ onto $N_0$.  Clearly $h$ maps
$\bar \psi$ to $\bar \psi$.  From now on we work in $V$.

We build an increasing continuous elementary chain $\langle N_\alpha:\alpha
< \omega_1 \rangle$, choosing $N_\alpha$ by induction on $\alpha$, as follows.
Note the $N_\alpha$'s are not necessarily
transitive or even well founded.

Let $\Gamma = \Gamma_\alpha = \{\varphi(v):N_\alpha \models 
``\{\delta \in h(S):\varphi(\delta)\}$ is stationary" and 
$\varphi \in \Phi_\alpha\}$ where $\Phi_\alpha$ is the set of first order 
formulas with parameters from $N_\alpha$ in the vocabulary $\{\in,<^*\}$
and the only free variable $v$.
Let $\le_{\Gamma_\alpha}$ be the following partial order of $\Gamma_\alpha:
\theta \le_{\Gamma_\alpha} \varphi$ iff $N_\alpha \models ``(\forall x)
[\varphi(x) \rightarrow \theta(x)]"$.  Let 
$t_\alpha$ be a subset of $\Gamma_\alpha$ such that:
\mr
\item "{$(a)$}"  $t_\alpha$ is downward closed, i.e. if $\theta
\le_{\Gamma_\alpha} \varphi$ and $\varphi \in t_\alpha$ \ub{then} 
$\theta \in t_\alpha$
\sn
\item "{$(b)$}"  $t_\alpha$ is directed
\sn
\item "{$(c)$}"  for some countable $M_\alpha \prec ({\Cal H}(\kappa),
\in,<^*)$ to which $N_\alpha$ belongs, if \nl
$\Gamma \in M_\alpha,\Gamma \subseteq \Gamma_\alpha$ is a 
dense subset of $\Gamma_\alpha$ \ub{then} 
$t_\alpha \cap \Gamma \ne \emptyset$.
\ermn
Clearly by the density if $\varphi \in \Gamma_\alpha$ and $\theta \in 
\Phi_\alpha$, then 
$\varphi \wedge \theta \in \Gamma_\alpha$ or $\varphi \wedge \neg \theta 
\in \Gamma_\alpha$.  Thus, $t_\alpha$ is a complete type over $N_\alpha$.
Since $N_\alpha$ has definable Skolem 
functions, we can let $N_{\alpha +1}$ be the Skolem hull of $N_\alpha \cup
\{b_\alpha\}$ where $N_\alpha \prec N_{\alpha +1},b_\alpha \in N_{\alpha +1}$
realizes $t_\alpha$.

We claim that $N_{\alpha +1}$ has no ``new natural numbers", i.e. if
$N_{\alpha +1} \models ``c$ is a natural numbers" then $c \in N_\alpha$.
Why?  As $c \in N_{\alpha +1}$ clearly for some $f \in N_\alpha$ we have
$N_\alpha \models ``f$ is a function with domain $\omega_1$, the countable
ordinals" and $N_{\alpha +1} \models ``f(b_\alpha) = c"$.  Let

$$
\align
{\Cal D}_f = \bigl\{ \varphi(v) \in \Gamma_\alpha:&N_\alpha 
\models ``(\forall x)(\varphi(x) \rightarrow f(x) 
\text{ is not a natural number})" \\
  &\text{or for some } d \in N_\alpha \text{ we have} \\
  &N_\alpha \models ``(\forall x)(\varphi(x) \rightarrow f(x) = d)" \bigr\}.
\endalign
$$
\mn
It is easy to check that ${\Cal D}_f$ is a subset of $\Gamma_\alpha$, it
belongs to $M_\alpha$ and it is a dense subset of $\Gamma_\alpha$; hence
$t_\alpha \cap {\Cal D}_f \ne \emptyset$.  Let $\varphi(x) \in {\Cal I}_f
\cap t_\alpha$, so $N_{\alpha +1} \models \varphi[b_\alpha]$, and by the
definition of ${\Cal I}_f$ we get the desired conclusion.

If $N_\alpha \models ``b$ is a countable ordinal" then $N_{\alpha +1} \models
``b < b_\alpha \and b_\alpha$ is a countable ordinal".  
Also $N_{\alpha +1} \models ``b_\alpha \in h(S)"$.

We claim that $b_\alpha$ is the least ordinal of $N_{\alpha +1} \backslash
N_\alpha$ in the sense of $N_{\alpha +1}$.  Assume
$N_{\alpha +1} \models ``c$ is a
countable ordinal, $c < b_\alpha"$ so for some $f \in N_\alpha$ we have
$N_\alpha \models ``f:\omega_1 \rightarrow \omega_1$ is a function" and
$N_{\alpha +1} \models ``c = f(b_\alpha)",N_{\alpha +1} \models ``f(b_\alpha) 
< b_\alpha"$.  
Then $N_\alpha \models ``\{\beta \in h(S):f(\beta) < \beta\}$ is a stationary
subset of $\omega_1$".  Let ${\Cal D} = 
\{\varphi(v) \in \Gamma_\alpha:(\exists \gamma < \omega_1)(\forall v)
(\varphi(v) \rightarrow f(v) = \gamma) \vee (\forall v)(\varphi(v) 
\rightarrow f(v) \ge v)\}$.  By Fodor's lemma (which $N_\alpha$
satisfies) ${\Cal D}$ is a dense subset of $\Gamma_\alpha$ and clearly 
${\Cal D} \in M_\alpha$.  Since $t_\alpha$ is sufficiently generic, there is
a $\gamma \in N_\alpha$ such that $N_{\alpha +1} \models ``f(b_\alpha) 
= \gamma"$.

Now $N_\alpha$ is not necessarily wellfounded but it has standard $\omega$
and without loss of generality
$N_\alpha \models ``a \subseteq \omega"$ implies $a = \{n < \omega:
N_\alpha \models ``n \in a"\}$ so as $h(\bar \psi) = \bar \psi$ 
clearly $N_\alpha \models ``x/E^{\bar \psi} \in B" \Rightarrow x/E^{\bar \psi}
\in B$, and $N_\alpha \models ``x,y,z \in B_*,x/E^{\bar \psi} + y/
E^{\bar \psi} = z/E^{\bar \psi}" \Rightarrow x/E^{\bar \psi}+y/E^{\bar \psi}
=z/E^{\bar \psi}$.

For each $\alpha < \omega_1$, if $N_\alpha \models ``b < \omega_1"$, let
$B^\alpha_b$ be the group $(h(\bar B))_b$ as interpreted in 
$N_\alpha$, i.e. $N_\alpha$ thinks that $B^\alpha_b$ is the $b$-th group in
the increasing chain $h(\bar B)$.  
Clearly $B^\alpha_b \subseteq B$ if $E^{\bar \psi}$ is the
equality, otherwise let $j^\alpha_b$ map $(x/E^{\bar \psi})^{N_\alpha}$ to
$x/E^{\bar \psi}$, so $j^\alpha_b$ embeds $B^\alpha_b$ into $B^0$; 
let this image be called $G^\alpha_b$.  Also in $N_\alpha$ there is 
a bijection between $B^\alpha_b$ and $\omega$.  If $\gamma > \alpha$, 
since $N_\alpha \preceq N_\gamma$ have the same natural numbers, clearly
$B^\alpha_b = B^\gamma_b$ when $E^{\bar \psi}$ is equality or 
$j^\alpha_b = j^\gamma_b$ and $G^\alpha_b = G^\gamma_b$ in the general case.  
In particular, $G^{\alpha +1}_{b_\alpha}$ is 
the union of $\{G^\alpha_b:N_\alpha \models ``b < \omega_1"\}$.

For $\alpha < \omega_1$, let $G_\alpha = G^{\alpha +1}_{b_\alpha}$ and
let $(h(\left< \langle b^\alpha_\ell:\ell \le n_\alpha \rangle:\alpha \in S
\right>))(b_\alpha) \in N_{\alpha +1}$ be $\langle (a^{b_\alpha}_\ell/
E^{\bar \psi})^{N_\alpha}:\ell \le m_\alpha \rangle$,
so $N_{\alpha +1}$ thinks that $\langle a^{b_\alpha}_\ell/E^{\bar \psi}:\ell
\le m_\alpha \rangle$ witness that $h(B)/B^{\alpha +1}_{b_\alpha}$ is not
free.  
Clearly $a^{b_\alpha}_0/E^{\bar \psi},\dotsc,a^{b_\alpha}_{m_\alpha}/
E^{\bar \psi} \in G_{\alpha +1}$ and

$$
PC(G_\alpha \cup \{a^{b_\alpha}_0/E^{\bar \psi},\dotsc,a^{b_\alpha}
_{m_\alpha}/E^{\bar\psi}\})/G_\alpha
$$
\mn
is not free.  So $G_{\alpha +1}/G_\alpha$ is not free.  Let $G = \dbcu
_{\alpha < \omega_1} G_\alpha$.  Then $G$ is not free.  But $G$ is a subgroup
of $B$, thus $B$ is not $\aleph_2$-free.  \hfill$\square_{\scite{1.2}}$
\enddemo
\bigskip

\remark{Remark}  Instead of the forcing we could directly build the 
$N_\alpha$'s but we have to deal with stationary subsets of ${}^\omega 2$ 
instead of $\omega_1$.
\endremark
\bigskip

\demo{\stag{1.3} Corollary}  If $B$ is an $\aleph_1$-free Borel 
abelian group, \ub{then} $B$ is $\aleph_2$-free if and only if 
$\{K \subseteq B:|K| =
\aleph_0$ and $B/K$ is $\aleph_1$-free$\}$ is not stationary.
\enddemo
\bn
\ub{\stag{1.4} Fact}:  If $2^{\aleph_0} < 2^{\aleph_1}$ then every Borel
Whitehead group $B$ is $\aleph_2$-free.
\bigskip

\demo{Proof}  By \cite{DvSh:65} (or see \cite{EM}) as
$2^{\aleph_0} < 2^{\aleph_1}$ we have: if $G$ be a 
Whitehead group of cardinality $\aleph_1$ and
$G = \dbcu_{\alpha < \omega_1} G_\alpha$ is such that 
$\langle G_\alpha:\alpha < \omega_1 \rangle$ is an increasing continuous 
chain of countable subgroups,
then $\{\alpha:G_{\alpha +1}/G_\alpha$ is not free$\}$ does not contain a
closed unbounded set (see \cite[Ch.XII,1.8]{EM}).  Thus, if $B$ is not
$\aleph_2$-free, then the subgroup $G$ constructed in the proof of lemma
\scite{1.2} is not Whitehead.  Since being Whitehead is a hereditary 
property (see \cite{EM}),
$B$ is not Whitehead.  \nl
${{}}$  \hfill$\square_{\scite{1.4}}$
\enddemo
\bn
The lemma shows that
\demo{\stag{1.5} Conclusion}   For Borel abelian groups
$B^{\bar \psi},``B^{\bar \psi}$ is $\aleph_2$-free" is absolute 
(in fact it is a $\sum^1_1$ property of $\bar \psi$).
\enddemo
\bigskip

\demo{Proof}  The formula will just say that there is a model of a suitable
fragment of ZFC (e.g. ZC) with standard $\omega$ to which $\bar \psi$ belongs
and it satisfies ``$B^{\bar \psi}$ is $\aleph_2$-free". \nl
${{}}$  \hfill$\square_{\scite{1.5}}$
\enddemo
\newpage

\head {\S2 On $\aleph_2$-free Whitehead} \endhead  \resetall
\bigskip

\proclaim{\stag{2.1} Theorem}  If $B$ is a Borel Whitehead group, \ub{then}
$B$ is $\aleph_2$-free.
\endproclaim
\bn
\ub{\stag{2.1A} Conclusion}: (CH) Every Whitehead Borel abelian group is free.
\bn
Before we prove we quote \cite[Definition 3.1]{Sh:44}.
\definition{\stag{2.2} Definition}  1) If $L$ is a subset of the 
$\aleph_1$-free abelian group, $G,PC(L,G)$ is the
smallest pure subgroup of $G$ which contains $L$.  Note that if $H$ is a pure
subgroup of $G,L \subseteq H$ then $PC(L,G) = PC(L,H)$.  We omit $G$ if it is
clear. \nl
2) If $H$ is a subgroup of $G,L$ a finite subset of $G,a \in G$, we say that
$\pi(a,L,H,G)$ means that: 
$PC(H \cup L) = PC(H) \oplus PC(L)$ but for no $b \in
PC(H \cup L \cup \{a\})$ is $PC(H \cup L \cup \{a\}) = PC(H) \oplus PC(L \cup
\{b\})$.
\enddefinition
\bigskip

\demo{Proof}  Assume $B$ is not $\aleph_2$-free.  We repeat the
proof of Lemma \scite{1.2}.  So in $V^{\bold P},B$ is a non-free 
$\aleph_1$-free abelian group of cardinality $\aleph_1$.  Hence by 
\cite[p.250,3.1(3)]{Sh:44}, $B$ satisfies possibility I or possibility II 
where we have chosen $\bar B = \langle B_\alpha:\alpha < \omega_1 \rangle$
increasing continuous with $B_\alpha$ countable, 
$B = \dbcu_{\alpha < \omega_1} B_\alpha$; the possibilities are
explained below.  The proof splits into the two cases.
\enddemo
\bn
\ub{Possibility I}:  By \cite[p.250]{Sh:44}.

So we can find (still in $V^{\bold P}$) an ordinal 
$\delta < \omega_1$ and $a^\ell_i
\in B$ for $i < \omega_1,\ell < n_i$ such that
\mr
\item "{$(A)$}"  $\{a^i_\ell + B_\delta:\ell < \omega_1,\ell \le n_i\}$ is
independent in $B/B_\delta$
\sn
\item "{$(B)$}"  $\pi(a^\ell_{n_i},L_i,B_\delta,B)$ where $L_i$ is the
subgroup of $B$ generated by $\{a^i_\ell:\ell < n_i\}$.
\ermn
This situation does not survive well under the process and the proof of Lemma
\scite{1.2} but after some analysis a revised version will.

Without loss of generality $n_i = n(*) = n^*$ (by the pigeon hole principle).
Let $N \prec ({\Cal H}(\chi),\in,<^*)$ be countable such that $\mu,
B_\delta,B,\langle B_\alpha:\alpha < \omega_1 \rangle,\left< \langle 
a^i_0,\dotsc,a^i_{n_i} \rangle:i < \omega_1 \right>$ belong to $N$.  
We can find $M \in V,M \cong N$; \wilog \, $M$ is transitive (so 
$M \models ``n$ is a natural number" iff $n$ is a natural number).

Let ${\frak B} \prec ({\Cal H}(\chi),\in,<^*)$ be countable, 
$M \in {\frak B}$.  Let $\Phi_M$ be the set of f.o. formulas $\varphi(v)$
in the vocabulary $\{\in,<^*\}$ and parameters from $M$ and the only free
variable $v$.  Now we imitate the proof of \cite{Sh:202}.  Let 
$\Gamma = \{\varphi(v) \in \Phi_M:M \models
``\{\alpha < \omega_1:\varphi(\alpha)\}$ is uncountable"$\}$
(equivalently $\Gamma$ is $\{a \subseteq \omega_1:|a| = \aleph_1\}^M$).
We can find $\langle t_\eta(v):\eta \in {}^\omega 2 \rangle$ such that:
\mr
\item "{$(a)$}"  each $t_\eta(v)$ a suitable generic subset of $\Gamma$,
i.e. $\Gamma$, is ordered by $\varphi_1(v) \le \varphi_2(v)$ if $M \models
(\forall v)(\varphi_2(v) \rightarrow \varphi_1(v))$ so $t_\eta(v)$ is
directed, downward closed and is not disjoint to any dense 
subset of $\Gamma$ from ${\frak B}$ \nl

\sn
\item "{$(b)$}"  for $k < \omega,\eta_0,\dotsc,\eta_{k-1} \in {}^\omega 2$
which are pairwise distinct \nl
$\langle t_{\eta_0}(v),\dotsc,t_{\eta_{k-1}}(v) \rangle$ 
is generic too (for $\Gamma^k$), i.e. if ${\Cal D} \in {\frak B}$ is a dense
subset of $\Gamma^k$ then $\dsize \prod_{\ell < k} t_{\eta_\ell}(v)$ is not
disjoint to ${\Cal D}$.
\ermn
(See explanation in the end of the proof of case II). \nl
So for each $\eta,t_\eta(v)$ is a complete type over $M$ hence we can find
$M_\eta,M \prec M_\eta,M_\eta$ the Skolem hull of $M \cup \{y_\eta\}$ such
that $y_\eta$ realizes $t_\eta(v)$ in $M_\eta$.  
So $M_\eta \models ``y_\eta$ a countable
ordinal".  Without loss of generality if $M_\eta \models ``\rho \in
{}^\omega 2"$ then $\rho \in {}^\omega 2$ and $\rho(n) = i \Leftrightarrow
M_\eta \models \rho(n)=i$ when $n < \omega,i < 2$.
\sn
Let $h:N \rightarrow M$ be the isomorphism from $N$ onto $M$.  We still use
$B_\delta$!  As $\bar a = \left< \langle a^i_\ell:\ell \le n^* \rangle:
i < \omega_1 \right> \in N$ we can look at $\bar a$ and $h(\bar a)$ as a
two-place function (with variables written as superscript and subscript).
So we can let $a^\eta_\ell(\ell \le n^*,\eta \in {}^\omega 2)$ be reals
such that: $M_\eta \models ``h(\bar a)^{y_\eta}_\ell = a^\eta_\ell"$. By
absoluteness $a^\ell_\eta \in B$ (more exactly $a^\ell_\eta \in B_* =
B^{\bar \psi}_*,a^\ell_n/E^{\bar \psi} \in B$) and 
$\pi(a^\eta_{n^*},\langle a^\eta_\ell:\ell < n^* \rangle,B_\delta,B)$.
\sn
If we can prove that $\langle a^\eta_\ell:\eta \in {}^\omega 2,\ell \le n^*
\rangle$ is independent over $B_\delta (= h(B_\delta))$, then the proof of
\cite[3.3]{Sh:98} finish our case: proving $B$ is not Whitehead group.  But
independence is just a demand on every finite subset.  So it is enough to
prove
\mr
\item "{$\otimes$}"  if $k < \omega,\eta_0,\dotsc,\eta_{k-1} \in
{}^\omega 2$ are distinct, then \nl
$\{a^{\eta_m}_\ell:\ell \le n^*,m < k\}$ is independent over $B_\delta$.
\ermn
We prove this by induction on $k$.  For $k=0$ this is vacuous, for 
$k=1$ it is part of the
properties of each $\langle a^\eta_\ell:\ell \le n^* \rangle$.  So let us
prove it for $k+1$.  Remember that $\langle t_{\eta_0}(v),\dotsc,
t_{\eta_k}(v) \rangle$ (more exactly $\dsize \prod_{\ell \le k}
t_{\eta_\ell}(v))$ is a generic subset of $\Gamma^k$.

Assume the desired conclusion fails.  So by absoluteness we can find
$\varphi_\ell(v) \in t_{\eta_\ell}(v)$ and $s^m_\ell \in \Bbb Z$ for $m \le k,
\ell \le n^*$ such that:
\mr
\item "{$\oplus$}"  if $t'_{\eta_m}(v) \subseteq \Gamma$ is generic over
${\frak B}$ for $m \le k$, moreover $\langle t'_{\eta_m}(v):m \le k \rangle$ 
is a generic subset of $\Gamma^k$ over ${\frak B}$ and $\varphi_m(v) \in 
t'_{\eta_m}(v)$, \ub{then} (defining $M'_{\eta_m}$ by $t'_{\eta_m}(v)$ and $a^{\eta_m}_\ell$
as before) 
$\dsize \sum \Sb \ell \le n^* \\ m \le k \endSb s^m_\ell a^{\eta_m}_\ell
= t \in B_\delta$.
\ermn
Clearly for $m \le k$ we have 
$M \models ``\{v:M \models ``\varphi_m(v) \wedge v \text{ a countable
ordinal}"\}$ has order type $\omega_1$" and \wilog \, $M \models ``\{v:M
\models ``\neg \varphi_m(v) \wedge v$ a countable ordinal"$\}$ has order type
$\omega_1"$.

So in $M$ there are $g_0,\dotsc,g_k \in M$ such that: $M \models
``g_i$ is a permutation of $\omega_1$, for $i \le k$ we have
$(\forall v)(\varphi_0(v) \leftrightarrow \varphi_0(g_i(v))$ and $g_0(v),
g_1(v),\dotsc,g_k(v)$ are pairwise distinct".  Let for $m \le k,
t^i_{\eta_0}(v) = \{\varphi(v) \in \Gamma:
\varphi(g_i(v)) \in t_{\eta_0}(v)\}$.  Let in $M_{\eta_0},
y^i_{\eta_0} = [g_i(y_{\eta_0})]^{M_{\eta_0}},a^{\eta_0,i}_\ell =
[h(\bar a)^{(y^i_{\eta_0})}_\ell]^{M_{\eta_0}}$.  Now $y^i_{\eta_0}$
realizes $t^i_{\eta_0}(v)$ and $M_{\eta_0}$ is also the Skolem hull of
$M \cup \{y^i_{\eta_0}\}$ and $\langle t^i_{\eta_0}(v),t_{\eta_1}(v),\dotsc,
t_{\eta_k}(v) \rangle \subseteq \Gamma^{k+1}$ is generic over
${\frak B}$ and $\varphi_0(v) \in t^i_{\eta_0}(v),\varphi_1(v) 
\in t_{\eta_1}(v),\dotsc,\varphi_k(v) \in t_{\eta_k}(v)$.  Hence for each 
$i \le k$ in $B$ we have $\dsize \sum_{\ell \le n^*} s^0_\ell  
a^{\eta_0,i}_\ell + \dsize \sum \Sb 0 < m \le k \\ \ell \le n^* \endSb 
s^m_\ell a^{\eta_m}_\ell = t \in B_\delta$.
\mn
By linear algebra $\{a^{\eta_0,i}_\ell:i \le k,\ell \le n^*\}$ is
not independent (actually, $i=0,1$ suffices - just subtract the 
equations).  By absoluteness this holds in $M_{\eta_0}$.  But the formula
saying this is false holds in $({\Cal H}(\chi),\in,<^*)$ hence in $N$, 
hence in $M$, hence in $M_\eta$ (it speaks on $\bar a,B,B_\delta$), 
contradiction.  So $\oplus$ fails hence $\otimes$ holds so we have 
finished Possibility I. 
\bn
\ub{Possibility II of \cite[p.250]{Sh:44}}:  In this case we have ``not
possibility I" but $S = \{\delta < \omega_1:\delta \text{ a limit ordinal
and there are } a^\delta_\ell$ for $\ell \le n_\delta \text{ such that }
\pi(a^\delta_{\eta_\delta},\langle a^\delta_\ell:
\ell < n_\delta \rangle_B,B_\delta,B)\}$ is stationary; all in $V^{\bold P}$.
Now \wilog \, we can find $\langle \alpha^\delta_n:n < \omega \rangle$ such 
that: $\alpha^\delta_n < \alpha^\delta_{n+1},\delta = \dbcu_{n < \omega}
\alpha^\delta_n$, and there are $y^\delta_m \in B_{\delta +1},t^\delta_m \in
B_{\alpha^\delta_n +1}$ and $s^\delta_{m,\ell} \in \Bbb Z$, 
(for $\ell < n_\delta$) such that:  
\mr
\item "{$\boxtimes(*)_0$}"  $y^\delta_0 = a^\delta_{n_\delta}$ and
\sn
\item "{$(*)_2$}"  $s^\delta_{m,n_\delta} 
y^\delta_{m+1} = \dsize \sum_{\ell < n^*}
s^\delta_{m,\ell} a^\delta_\ell + y^\delta_m + t^\delta_m$ 
\sn
\item "{$(*)_3$}"  $s^\delta_{m,n_\delta} > 1$, morever if $s$ is a proper
divisor of $s^\delta_{m,n_\delta}$ (e.g. 1) \ub{then} 
$sy^\delta_{m+1,n_\delta}$ is not in $B_\delta + \langle 
\{a^\delta_i:\ell < n_\delta\} \cup \{y^\delta_m\} \rangle_B$
\sn
\item "{$(*)_4$}"  if $\alpha \in \delta \backslash \{\alpha^\delta_n:n <
\omega\}$ \ub{then} $PC_B(B_{\alpha +1} \cup \{a^\delta_0,\dotsc,
a^\delta_{n_\delta}\}) =$ \nl
$PC_B(B_\alpha \cup \{a^\delta_0,\dotsc,
a^\delta_{n_\delta}\}) + B_{\alpha +1}$ \nl
\sn
[why?  known, or see later.]
\ermn
Without loss of generality $\delta \in S \Rightarrow n_\delta = n^*$.
So as in the proof of Lemma \scite{1.2} we can choose countable
$N \prec ({\Cal H}(\chi),\in,<^*)$
such that $\bar a = \left < \langle a^\delta_\ell:\ell \le n^* \rangle:
\delta \in S \right>,\bar \alpha = \left < \langle \alpha^\delta_n:n < \omega
\rangle:\delta \in S \right>,\left< (\langle s^\delta_{m,\ell}:\ell \le
n^* \rangle,y^\delta_m,t^\delta_m)_{m < \omega}:\delta \in S \right>$ 
belongs to $N$, then define $M$ and choose ${\frak B}$ as before.  
We let this time $\Gamma = \Gamma_M$ be as in 
the proof of Lemma \scite{1.2}, that is $\{\varphi(v):M \models 
``\{\delta \in S:\varphi(\delta)\} \text{ stationary}\}$. \nl
We can find $\langle t_\eta(v):\eta \in {}^\omega 2 \rangle$ such that:
\mr
\item "{$(a)$}"  each $t_\eta(v) \subseteq \Gamma$ is generic over 
${\frak B}$ as before hence 
\sn
\item "{$(b)$}"  for $k < \omega$ and pairwise distinct
$\eta_0,\dotsc,\eta_{k-1} \in {}^\omega 2,\langle t_{\eta_0},\dotsc,
t_{\eta_{k-1}} \rangle$ is generic over ${\frak B}$
\sn
\item "{$(c)$}"  letting $M_\eta,y_\eta$ be such that: $M \prec M_\eta,
M_\eta$ the Skolem hull of $M_\eta \cup \{y_\eta\},y_\eta$ realizes
$t_\eta(v)$ in $M_\eta$ we have
\sn
{\roster
\itemitem{ $(i)$ }  $M_\eta \models ``y_\eta$ is a countable ordinal $\in
S"$
\sn
\itemitem{ $(ii)$ }  $M \models ``a$ is a countable ordinal" $\Rightarrow
M_\eta \models ``a < y_\eta"$
\sn
\itemitem{ $(iii)$ }  if $y \in M_\eta$ satisfies (i) + (ii) then $M_\eta
\models ``y_\eta < y"$.
\endroster}
\ermn
So looking at $h:N \rightarrow M$ the isomorphism, \ub{then} $\alpha^\eta_n =:
[h(\bar \alpha)]^{y_\eta}_n$ for $n < \omega$ satisfies:

$$
M_\eta \models ``\alpha^\eta_n \text{ a countable ordinal}"
$$

$$
M_\eta \models ``\alpha^\eta_n < \alpha^\eta_{n+1} < y_\eta"
$$

$$
M_\eta \models ``[h(\bar \alpha)]^y_\eta \text{ is unbounded below } y_\eta"
$$
\mn
hence $\{\alpha^\eta_n:n < \omega\} \subseteq M$ is unbounded among 
the countable ordinals of $M$. \nl
Now by easy manipulation (see proof below):
\mr
\item "{$(c)$}"  if $\eta_1 \ne \eta_2$ then $\{\alpha^{\eta_1}_n:n <
\omega\} \cap \{\alpha^{\eta_2}_n:n < \omega\}$ is finite.
\ermn
(We can be lazy here demanding just that no $\{\alpha^\eta_n:n < \omega\}$
is included in the union of a finite set with the union of 
finitely many sets of the
form $\{\alpha^\nu_n:n < \omega\}$ which follows from pairwise generic, 
and one has to do slightly more abelian group theory work below).
\mn
Now we can let $a^\eta_\ell = [(h(\bar a))^{y_\eta}_\ell]^{M_\eta}$.  By
linear algebra we get the independence hence a contradiction to our being
in possibility II (or directly get $\otimes$ in the proof in the case
possibility I holds). \nl
An alternative is the following:
\sn
We are assuming that in $V^{\bold P}$, possibility I fails.  So also in $V$,
letting $A = M \cap B^{\bar \psi}$ the following set is countable: $K[A] =:
\{\langle a_\ell:\ell \le n \rangle:n < \omega,a_\ell \in B,\langle a_\ell:
\ell \le n \rangle \text{ independent over } A \text{ in } B \text{ and }
\pi(a_n,\langle a_\ell:\ell < n \rangle_B,A,B)\}$ (see proof later). \nl
For each such $\bar a = \langle a_\ell:\ell \le n \rangle$ we can look at a
relevant type it realizes over $A$

$$
\align
t(\bar a,A) = \bigl\{(\exists y)(sy = \dsize \sum_{\ell \le n} s_\ell
x_\ell):&B \models (\exists y)(sy = \sum s_\ell a_\ell), \\
  &s,s_\ell \text{ integers} \bigr\}
\endalign
$$
\mn
so $\{t(\bar a,A):\bar a \in K[A]\}$ is countable.  But for the $\eta \in
{}^\omega 2$ the types \nl
$t(\langle a^\eta_\ell:\ell < n_\eta \rangle,A)$ are pairwise 
distinct, contradiction, so actually case II never occurs.
\bn
We still have some debts in the treatment of possibility II. \nl
\ub{Why do clauses (b) and (c) hold}?  For each $n$ we let

$$
\align
\Gamma_{M,n} = \biggl\{ \varphi(v):&(i) \quad \varphi(v) \text{ is a first
order formula with parameters from } M \\
  &(ii) \quad \text{for some } \beta^*_\ell \in M \cap \omega_1 \text{ for }
\ell < n \text{ we have} \\
  &\qquad \quad M \models ``(\forall v)(\varphi(v) \rightarrow v \in h(S))
\and \dsize \bigwedge_{\ell < n} (h(\bar \alpha))^v_\ell = \beta^*_\ell) \\
  &(iii) \quad M \models ``(\forall \beta < \omega_1)
(\exists^{\text{stat}} v < \aleph_1)
[(\varphi(v) \and \beta < (h(\bar \alpha))^v_n)]" \biggr\}.
\endalign
$$
\mn
Now note:
\mr
\item "{$\otimes_0$}"  $\Gamma_{M,n} \subseteq \Gamma_M$
\sn
\item "{$\otimes_1$}"  if $\varphi(v) \in \Gamma_M$ and $n < \omega$ then
for some $m \in [n,\omega)$ and $\beta_\ell \in M \cap \omega_1$ for $\ell
< m$ we have $``\varphi(v) \and 
\dsize \bigwedge_{\ell < m} ``(h(\bar \alpha))^v_\ell = \beta_\ell"$ 
belongs to $\Gamma_{M,m}$
\sn
\item "{$\otimes_2$}"  if $\varphi(v) \in \Gamma_{M,n}$ and $\beta \in M \cap
\omega_1$ then $\varphi'(v) = \varphi(v) \and \beta < (h(\bar \alpha))^v_n$
belongs to $\Gamma_{M,n}$.
\ermn
Now let $\langle {\Cal D}_n:n < \omega \rangle$ be the family of 
dense open subsets
of $\Gamma_M$ which belong to ${\frak B}$.  We choose by induction on $n,
\langle \varphi_\eta(v):\eta \in {}^n 2 \rangle,k_\eta < \omega$ such that:
\mr
\item "{$(\alpha)$}"  $\varphi_n(v) \in \Gamma_{M,k_\eta}$
\sn
\item "{$(\beta)$}"  $\varphi_\eta(v) \in {\Cal D}_\ell$ if 
$\ell < \ell g(\eta)$
\sn
\item "{$(\gamma)$}"  $\varphi_\eta(v) \le_\Gamma \varphi_{\eta \char 94
\langle i \rangle}(v)$ for $i=0,1$
\sn
\item "{$(\delta)$}"  if $\eta_0 \ne \eta_1 \in {}^n 2,\eta_i \triangleleft
\nu_i \in {}^{n+1}2$ for $i=0,1$ and 
$k_{\eta_0} \le k < k_{\nu_0}$ and $M \models
(\forall v)(\varphi_{\nu_0}(v) \rightarrow (h(\bar \alpha))^v_k = \beta)$
\ub{then} $M \models (\forall v)[\varphi_{\nu_1}(v) \rightarrow
\dsize \bigwedge_{\ell < k_{\nu_1}} (h(\bar \alpha))^v_\ell \ne \beta]$.
\ermn
There is no problem to do it and 
$t_\eta(v) = \{\varphi(v) \in \Gamma_M:\varphi(v)
\le_{\Gamma_M} \varphi_{\eta \restriction n}(v)$ for some $n < \omega\}$ for
$\eta \in {}^\omega 2$ are as required.
\bn
\ub{Why does $\boxtimes$ hold}?

For $\delta \in S$ let $w_\delta = \{\alpha < \delta:PC_B(B_{\alpha +1} \cup
\{a^\delta_0,\dotsc,a^\delta_{n,\alpha}\})$ is not equal to
$PC_B(B_\alpha \cup \{a^\delta_0,\dotsc,a^\delta_{n,\alpha}\}) +
B_{\alpha +1} \subseteq B\}$.

Let $S' = \{\delta \in S:(\forall \alpha < \delta)(|w_\delta \cap \alpha|
< \aleph_0)\}$, if $S'$ is stationary we get $\boxtimes$, otherwise $S
\backslash S'$ is stationary, and for $\delta \in S \backslash S'$ let
$\alpha_\delta = \text{ Min}\{\alpha:w_\delta \cap \alpha$ is infinite$\}$.
By Fodor's lemma for some 
$\alpha(*) < \omega_1,S'' = \{\delta \in S \backslash
S':\alpha_\delta = \alpha(*)\}$ is stationary hence uncountable and we can
get possibility I, contradiction.  \hfill$\square_{\scite{2.1}}$
\newpage

\head {\S3 Refinements} \endhead  \resetall
\bn
We may wonder if we can weaken the demand ``Borel".
\definition{\stag{3.3} Definition}
1) We say $\bar \psi$ is a code for a
Souslin abelian group if in Definition \scite{0.1} we weaken the demand on
$\psi_0,\psi_1$ to being a $\sum^1_1$ relation. \nl
2) A model $M$ of a fragment of ZFC is essentially transitive if:
\mr
\item "{$(a)$}"  if $M \models ``x$ is an ordinal" and $(\{y:y <^M x\},
\in^M)$ is well ordered then $x$ is an ordinal and $M \models ``y \in x"
\Leftrightarrow y \in x$
\sn
\item "{$(b)$}"  if $\alpha$ is an ordinal, $(\{y:y <^M x\},\in^M)$ is well
ordered and $M \models ``\alpha$ an ordinal,
rk$(x) = \alpha"$, then $M \models ``y \in x" \Leftrightarrow y \in x$.
\ermn
3) For $M$ essentially transitive with standard $\omega$ such that
$\bar \psi \in M$ let $B^M$ is $B^{\bar \psi}$ as interpreted in $M$ and
trans$(M) = \{x \in M:x$ as in (b) of part (2)$\}$.

\enddefinition
\bigskip

\demo{\stag{3.4} Fact}  1) ``$\bar \psi$ codes a Souslin abelian group" in a
$\Pi^1_2$ property. \nl
2) If $M$ is a model of a suitable fragment of set theory (comprehension is
enough), \ub{then} $M$ is isomorphic to an essentially transitive model. \nl
3) If $M$ is an essentially transitive model with standard $\omega$ of a
suitable fragment of ZFC and $\bar \psi \in M$, (note $\bar \psi$ is really
a pair of subsets of ${\Cal H}(\aleph_0))$, \ub{then} letting
$B^{\bar \psi} = (B^{\bar \psi})^M \cap$ trans$(M)$ there is a homomorphism
$\bold j_M$ from $B^M$ into $B = B^{\bar \psi}$ such that $M \models ``t =
x/E^{\bar \psi}"$ implies $\bold j_M(t) = x/E^{\bar \psi}$. \nl
4) If $M \prec N$ are as in (3), then $\bold j_M \subseteq \bold j_N$.
\enddemo
\bigskip

\demo{Proof} Straightforward.
\enddemo
\bigskip

\proclaim{\stag{3.5} Claim}  1) In \scite{1.2}, \scite{2.1} we can assume 
that $B = B^{\bar \psi}$ is only Souslin. \nl
2) If $B = B^{\bar \psi}$ is not $\aleph_2$-free, \ub{then} case I of 
\cite{Sh:44}(3.1) holds, more of the conclusion of case I in the proof of
\scite{2.1} holds.
\endproclaim
\bigskip

\remark{Remark}  If only $\psi_1$ is Souslin, i.e. is $\sum^1_1$, just repeat
the proofs.
\endremark
\bigskip

\demo{Proof}  For both we imitate the proof of \scite{2.1}.

In both possibilities, for each $\eta \in {}^\omega 2$, let $G_\eta$ be
the group which $\bar \psi$ defines in $M_\eta$, (the $M_\eta$'s chosen as
there).  So $\bold j_{M_\eta}$ is a homomorphism from $G_\eta$ into $B$.
However, $\bold j_M \subseteq \bold j_{M_\eta}$ and $\bold j_M$ is one to
one.  Now in defining $\pi(x,L,B_\delta,B)$ we can add that we cannot find
$L' \cup \{x'\} \subseteq PC(B_\delta \cup L \cup \{x\})$ such that
$\pi(x',L',B_\delta,B)$ and $|L'| < |L|$, i.e. the $n$ is minimal.  As $B$
is $\aleph_1$-free, this implies that $\bold j_M \restriction B(PC(B_\delta
\cup \{a^n_\ell:\ell \le n^*\})^{M_\eta}$ is one to one and by easy algebraic
argument, we can get, for \scite{2.1}, non-Whiteheadness and for \scite{1.2},
non $\aleph_2$-freeness.  \hfill$\square_{\scite{3.5}}$
\enddemo
\bigskip

\demo{\stag{3.6} Fact}  1) ``$B^{\bar \psi}$ is non-$\aleph_2$-free" is a
$\sum^1_1$-property of $\bar \psi$, assuming $B^{\bar \psi}$ is a
$\aleph_1$-free Souslin abelian group. \nl
2) ``$\bar \psi$ codes a $\aleph_1$-free Souslin abelian group" is a
$\Pi^1_2$-property of $\bar \psi$.
\enddemo
\bigskip

\demo{Proof}  Just check.
\enddemo
\newpage

\shlhetal

\newpage
    
REFERENCES.  
\bibliographystyle{lit-plain}
\bibliography{lista,listb,listx,listf,liste}

\enddocument

\bye